\documentclass{amsart}
\usepackage{amsmath, amsfonts,amssymb, latexsym}

\newtheorem{theorem}{Theorem}[section]
\newtheorem{cor}[theorem]{Corollary}

\newtheorem{thm}[theorem]{Theorem}
\newtheorem{definition}[theorem]{Definition}

\newtheorem{remark}[theorem]{Remark}
\newtheorem{rmk}[theorem]{Remark}
\newtheorem{lemma}[theorem]{Lemma}

\newtheorem{prop}[theorem]{Proposition}

\newcommand{\ag}{{\alpha}}
\newcommand{\al}{{\alpha}}

\newcommand{\Og}{{\Omega}}

\newcommand{\om}{{\omega}}

\newcommand{\io}{{\iota}}

\newcommand{\ga}{{\gamma}}

\newcommand{\la}{{\lambda}}
\newcommand{\La}{{\Lambda}}
\newcommand{\si}{{\sigma}}
\newcommand{\tg}{{\theta}}

\newcommand{\1}{{{\mathchoice {\rm 1\mskip-4mu l} {\rm 1\mskip-4mu l}{\rm 1\mskip-4.5mu l} {\rm 1\mskip-5mu l}}}}

\newcommand{\GW}[3][M]{ \ensuremath{{\rm GW}^{#1}_{#2,#3}}}
\newcommand{\bb}[1]{\ensuremath{\mathbb{#1}}}

\newcommand{\QED}{\hfill$\Box$\medskip}
\newcommand{\mo}[2]{\overline{\mathcal{M}}_{0,2}(\widetilde{M}_\la, \tilde{J}, \si_I + B;#1,#2)}
\newcommand{\vmc}[1][3]{{\overline{\mathcal{M}}_{0,#1}^\nu}}
\newcommand{\moduli}[1][3]{\overline{\mathcal{M}}_{0,#1}}
\newcommand{\modu}[2][0]{{\mathcal{M}}_{#1,#2}}
\newcommand{\T}[1]{\ensuremath{\widetilde{#1}}}
\renewcommand{\to}{\longrightarrow}
\newcommand{\inject}{\hookrightarrow}
\newcommand{\PD}{{\rm PD} }
\newcommand{\eff}{\ensuremath{K^{\rm eff}}}
\newcommand{\x}[1]{\ensuremath{x_{#1}}}
\newcommand{\xc}[1]{\ensuremath{x_{#1^c}}}
\newcommand{\ee}[1]{\ensuremath{\otimes e^{A_{#1}}}}
\newcommand{\e}[1]{\ensuremath{\otimes e^{#1}}}
\newcommand{\sei}[1]{\Psi(\lambda,{#1})}
\newcommand{\xp}{x_{\mathcal{S}}}
\newcommand{\p}[2]{p^{#1}_{#2}}
\newcommand{\grad}{\ensuremath{{\rm grad}}}

\newcommand{\simax}{\si_{\rm max}}

\begin{document}
  
  \title{Quantum Cohomology and $S^1$-actions with isolated fixed points}
  \author{Eduardo Gonzalez}
    \thanks{partial support by CONACyT-119141}
  \email{eduardo@math.sunysb.edu}
  \address{Mathematics Department\\
    SUNY at Stony Brook\\
    Stony Brook NY, 11777
  }
  \date{October, 2003}
  \begin{abstract}
    This paper studies symplectic  manifolds that admit semi-free circle
    actions with isolated fixed points. We prove, using results on the Seidel 
    element \cite{McT}, that the (small)
    quantum cohomology of a $2n$ dimensional manifold of this type is
    isomorphic  to the (small) quantum cohomology of a product of $n$ 
    copies of $\bb{P}^1$. This generalizes a result due to Tolman and Witsman \cite{TW}. 
  \end{abstract}
  \maketitle
  \section{Introduction}\label{intro}
  
  Let $(M,\om)$ be a $2n$ dimensional compact, connected, symplectic manifold, and 
  let $\{\la_t\}=\la: S^1
  \to \mathrm{Symp}(M,\om)$ be a symplectic circle action on $M$, that is, if $X$
  is the vector field generating the action, then $\mathcal{L}_X \om = d \io_X
  \om = 0$. Recall that the action is semi-free if it is free on $M \backslash M^{S^1}$. 
  This is equivalent to say that the only {\em weights} at every fixed point are 
  $\pm 1$. A circle action is said to be  Hamiltonian if  there is  
  a $C^{\infty}$  function $H:M \to \bb{R} $ such that $\io_X \om = -d H$. Such a function is 
  called a Hamiltonian  for the action.
  
  Tolman and Weitsman proved in \cite{TW} that if the action is semi-free and admits only 
  isolated fixed points, then  the action must be Hamiltonian provided that there is at
  least one fixed point. 
  There is a great deal of information concerning the topology of manifolds carrying such 
  actions. The first result in this direction is due 
  to Hattori~\cite{Ha}. He proves that there is an isomorphism from the cohomology ring
  $H^*(M;\bb{Z})$ to the cohomology ring of a product of $n$ copies of $\bb{P}^1$.
  Moreover, this isomorphism preserves Chern classes. In  \cite{TW} Tolman and Weitsman 
  generalize Hattori's result  to equivariant cohomology. 
  The main result of this paper is to extend this result to  quantum cohomology. 
  In \S\ref{ss:qh} we prove that $M$ is almost Fano manifold, therefore we can
  use polynomial coefficients $\La:=\bb{Q}[q_1,\dots,q_n]$ for the quantum cohomology
  ring. The main theorem is the following.

  \begin{thm}\label{main} Let $(M,\om)$ be a $2n$-dimensional compact connected symplectic 
    manifold. Assume $M$ admits a semi-free circle action with a finite non-empty set of
    fixed points. Then there is an isomorphism of (small) quantum cohomology 
    $$QH^*(M;\La)\cong QH^*((\bb{P}^1)^n;\La).$$
  \end{thm}
   
  Note that we can compute directly the quantum cohomology  of  $\bb{P}^1 \times \dots 
  \times \bb{P}^1$ to get the following result. 

  \begin{cor}\label{toric} The (small) quantum cohomology of $M$ is given by 
    $$
    QH^*(M;\La)\cong QH^*((\bb{P}^1)^n;\La)\cong \frac{\bb{Q}[x_1 ,\dots,
    x_n,q_1,\dots,q_n]}{ <x_i*x_i-q_i>} 
    $$
    where $\deg (x_i)=2$ and $\deg q_i=4$.
    
    Moreover, all other products are given by   
    $$x_{i_1}* \dots *x_{i_k}=x_{i_1}\smile \dots \smile x_{i_k}$$ 
    for $i_1< \dots < i_k$. Here the product on the left is the quantum product, while the 
    term on the right is the usual cup product.
  \end{cor}

  To prove Theorem \ref{main} we will construct a set  of generators $\{ x_i\}$ of the
  cohomology ring 
  $H^*(M;\bb{Z})$. Then we prove in Lemma \ref{l:main} that the quantum products of these
  generators satisfy the expected relations given in Corollary \ref{toric}. 

  To get this
  relations we use a result of McDuff-Tolman \cite{McT} to understand how  the Seidel
  automorphism  acts on the generators. We will see in Corollary \ref{seidelact} that this
  action do not have higher order terms, that is the automorphism is given by single
  homogeneous terms in quantum cohomology. Thus the Seidel
  automorphism acts by  permutation of the elements in the basis.
  To construct such generators for the cohomology ring, we will adopt the tools that
  Tolman and Weitsman developed to prove the following theorem.

  \begin{thm}[\cite{TW}]\label{thm:TW} Let $(M,\om)$ be a compact, connected symplectic
    manifold with a semi-free, Hamiltonian circle action with isolated fixed points.
    Then, there is an isomorphism of rings $H^*_{S^1}(M) \simeq
    H^*_{S^1}(((\bb{P}^1)^n))$ which takes the equivariant Chern classes of 
    $M$ to those of $(\bb{P}^1)^n$. Therefore the equivariant cohomology ring is given by 
    
    $$H^*_{S^1}(M)= \bb{Z}[a_1 ,\dots, a_n, y]/ (a_i y - a_i^2).$$ 
    Here $a_i\in H^2_{S^1}(M)$  and the
    equivariant Chern series is given by $c_t(M)=\sum_i{c_i(M) t^i}$ where 
    
    $$
    c_t(M)=\prod_i(1+t(2a_i - y)).
    $$ 
    
  \end{thm}

  Although Tolman and Weitsman use equivariant cohomology for getting an invariant base
  for  $H^*(M;\bb{Z})$, the results of McDuff-Tolman require a more geometric description
  of the basis. Therefore the crucial element in most of the results of this  paper is
  having geometric representatives of the cycles dual to the cohomology basis. These
  geometric representatives are defined by the Morse complex of the Hamiltonian function.    

  The paper is organized as follows. All the Morse theoretical constructions are in
  \S\ref{ss:morse}. In section \ref{ss:eqc} we use equivariant cohomology to provide an
  invariant basis for cohomology. Then we establish the relation with the Morse cycles. In
  \S\ref{ss:qh} we define the quantum cohomology ring and we get results that help to
  reduce the quantum product formulas. In \S\ref{ss:seidel} we define the Seidel
  automorphism in quantum cohomology. In \S\ref{ss:seidel2} we relate the Seidel
  automorphism with invariant chains. Then we compute explicitly the Seidel
  element. Finally in \S \ref{ss:proof} we use the associativity of the quantum product
  together with some dimensional arguments to provide the proof of Theorem \ref{main}. 
  
  \medskip
  
  \noindent {\bf Acknowledgments:} The author thanks Dusa McDuff for all her
  encouragement, patience, generosity and support. It would be impossible to finish this
  work without her help. The author also  thanks CONACyT for their support.

  \medskip
  
  \section{Morse Theory and  Equivariant Cohomology}\label{s:tools}
  In this section we establish all the tools we need to prove Theorem \ref{main}. We start
  in \S\ref{ss:morse} with basic definitions of Morse theory. For more
  details the reader can consult \cite{AuB,Sch} for Morse theory.
 
  Following the approach of  \cite{McT}, we will construct invariant
  Morse cycles to be able to calculate the Seidel element of $M$. This will be done in the
  next section. We introduce equivariant cohomology to identify a basis in cohomology and describe
  the relation with Morse cycles. At the end, we provide several results that will be
  necessary in \S \ref{ss:proof}.

  \medskip
  \subsection{Morse Theory}\label{ss:morse}
  As in \S \ref{intro}, let  $(M,\om)$ be a symplectic $2n$-dimensional manifold with
  $S^1$ action generated by a Hamiltonian function  $H$. Thus ~$\io_X\om=-d 
  H$ and  $X=J\text{grad}(H)$, where the gradient is taken respect to the metric
  $g_J(x,y)=\om(x,Jy)$ for an $\om$-compatible $S^1$-invariant almost complex structure
  $J$. With respect to this metric, $H$ is a (perfect) Morse function \cite{K} and the zeroes
  of $X$ are exactly the critical points of $H$.
  For each fixed point $p\in M^{S^1}$, denote by $\ag(p)$ the index of $p$ and let $m(p)$ be
  the sum of weights at $p$. Since the action is semi-free $m(p)=n_+(p)- n_-(p)$ where
  $n_+(p)$ is the number of positive weights and $n_-(p)$ the number of negative ones.
  Then  $\ag(p)=2n_-(p)=n-m(p)$.
  
  In order to understand the (co)homology of $M$ in terms of $S^1$-invariant cycles,
  we will consider the {\em stable} and {\em unstable} manifolds with respect to the
  gradient 
  flow $-\grad(H)$. More precisely, let $p,q$ be  critical points of $H$. Define the stable
  and unstable manifolds  by 
  
  $$W^s(q)=\{ \ga :\bb{R} \to M | \lim_{t \to \infty} \ga(t) = q\},$$ 
  $$W^u(p)=\{ \ga :\bb{R} \to M | \lim_{t \to - \infty} \ga(t) = p\}.$$
  Here $\ga(t)$ satisfies   the gradient flow equation 
  $$
  \ga'(t)=-\grad H(\ga(t)).
  $$
  These spaces are manifolds of dimension 
  $$
  \dim W^s(q)=2n-\ag(q) \text{ \ \ and \ \  } \dim W^u(p)=\ag(p),
  $$
  and the evaluation map  $\ga \mapsto \ga(0)$ induces smooth embeddings into $M$
  $$
  E_q:W^s(q) \to M \text{\  \ and \ \ }E_p:W^u(p) \to M. 
  $$
  
  When these manifolds intersect transversally for all fixed points $p,q$, the gradient
  flow is 
  said to be Morse-Smale \cite{AuB,Sch}. Under this circumstance we say that the pair
  $(H,g_J)$ is {\em Morse regular}.

  In \cite{Sch} Schwartz 
  proved that there is a way of  {\em partially compactifying} these manifolds and that there
  are natural extensions  of the evaluation maps so that 
  these compactifications with their evaluation maps $E_p: \overline{W^s(p)}\to M$ and
  $E_q:\overline{W^u(q)}\to M$, 
  define  {\em pseudocycles}.  The
  compactification of $W^s(p)$ is made by adding {\em broken trajectories} through fixed
  points of index $\ag(p)-1$.  When the action is semi-free and
  admits isolated fixed points, all the fixed points have even index (see comment after
  Theorem \ref{th:K}), therefore  $W^s(p)$ is already  compact in the sense of
  Schwartz. Thus $W^s(p)$ is itself a pseudocycle. The same is true for $W^s(x)$. It is
  well known that pseudocycles  define classes in homology (see \cite{McS2}).  We will
  denote by 
  $[W^u(p)]\in H_{\ag(x)}(M;\bb{Z})$ and  $[W^s(p)]\in H_{n-\ag(x)}(M;\bb{Z})$ the
  homology classes defined by these manifolds. To make these
  classes $S^1$-invariant we need to consider a special type of almost 
  complex structure, as we explain below. 

  Assume $(M, \om)$ admits a Hamiltonian $S^1$-action with isolated fixed points. Each
  fixed point $p\in M$ has a neighborhood $U(p)$ that is diffeomorphic to a 
  neighborhood of zero  in  a $2n$-dimensional Hermitian vector space $E(p)=E_1 \oplus
  \dots \oplus E_n$, in such 
  a way that the moment map $H$ is given by

  $$
  H(v_1, \dots , v_n)= \sum_j {\pi m_j |v_j|^2}
  $$
  and $S^1$ acts in $E_j$ just by multiplication by $e^{2\pi i m_j}$. Here the numbers
  $m_j\in \bb{Z}$ are exactly the weights of the action. Under the identification above,
  the almost-complex structure $J$  is the standard complex structure on the Hermitian
  vector space $E(p)$. Observe that $E(p)$ can be
  written as $E^+ \oplus E^-$ where $E^{\pm}$ is the sum of the $E_j$ where $m_j>0$ or
  $m_j<0$ respectively. We can call the spaces $E^{\pm}$  the positive and negative
  normal bundles to the point $p$.

  If we start with any compatible almost complex structure $J$ near the fixed
  points, we can extend $J$ to an $S^1$-invariant $\om$-compatible almost complex
  structure $J_M$ on $M$ whose restriction to the open sets $U(p)$ is $J$. Denote by
  $\mathcal{J}_{\rm inv}(M)$ the set of all $J$ that are equal to $J_M$ near the fixed
  points. 

  The following lemma shows that it is possible to acquire regularity with generic
  almost-complex structures.

  \begin{lemma}[\cite{McT}] Suppose that $H$ generates a semi free $S^1$-action on
  $(M,\om)$. Then for a generic choice of $J\in\mathcal{J}_{\rm inv}(M)$ the pair
  $(H,g_J)$ is Morse regular. 
  \end{lemma}

  For the rest of this paper, we will only consider Morse regular  pairs $(H,g_J)$ as in
  the previous lemma.

  \medskip
  
  \subsection{Equivariant Cohomology}\label{ss:eqc}
  We can start with a quick review of equivariant cohomology. Let $ES^1$ be  a contractible
  space where $S^1$ acts freely, and denote $BS^1=ES^1/S^1$. 
  Then $H^*(BS^1;\bb{Z})$ is the polynomial ring $\bb{Z}[y]$ where $y\in
  H^2(BS^1;\bb{Z})$. 

  Let $S^1$ act  on a manifold $M$. The equivariant cohomology of $M$, denoted by
  $H^*_{S^1}(M)$ is 
  defined by $H^*(M\times_{S^1} ES^1;\bb{Z})$. Note that $H^*(BS^1;\bb{Z})$ is naturally
  isomorphic to $H^*_{S^1}(pt)$, if $pt\in M$ is a point. Under this construction, we
  have two 
  natural maps, the projection $p:M\times_{S^1} ES^1 \to BS^1$ and the inclusion (as
  fiber) $i:M \to M\times_{S^1} ES^1$. The pullback $p^*:H^*(BS^1;\bb{Z})\to
  H^*_{S^1}(M)$ makes $H^*_{S^1}(M)$ a $H^*(BS^1;\bb{Z})$ module, while the restriction
  $i^*:H^*_{S^1}(M) \to H^*(M)$ is the ``reduction'' of invariant data to ordinary
  data. An immediate consequence is that $i^*(y)=0$. 
  
  Let $j:M^{S^1}\to M$ be the natural inclusion. In \cite{K}  Kirwan proved that if the
  action is Hamiltonian, the induced map  $j^*:H^*_{S^1}(M) \to  H^*_{S^1}(M^{S^1})$ is
  injective. The proof of this theorem is based on the following result, where we weaken
  the statement to  match our needs. For a fixed point $p\in M^{S^1}$ we denote by
  $a|_p:=(j_p)^*(a)$ where $(j_p)^*:H^*_{S^1}(M) \to H^*_{S^1}(p)$ and
  $j_p$ is the obvious inclusion.   
  
  \begin{thm}[\cite{K}]\label{th:K}
  Let  the circle act on a symplectic manifold $M$ in a Hamiltonian way. Assume the action is
  semi-free and that there are only isolated fixed points. Let $p\in M$ be a fixed point
  of index $2k$. Then there exists a unique class 
  $a_p\in H^{2k}_{S^1}(M)$ such that $a_p|_p=(-1)^k y^k$, and $a_p|_{p'}=0$ for all other fixed
  points $p'$ of index less than or equal to $2k$. Moreover, if we consider all fixed
  points, the  classes $a_p$ form a basis for $H^*_{S^1}(M)$ as a $ H^*(BS^1;\bb{Z})$
  module.  
  \end{thm}
  As a remark on the previous theorem, note that the term $(-1)^k y^k$ is the equivariant
  Euler class of the negative  normal bundle at $p$.
  
  As stated in \S \ref{intro},  there is
  an isomorphism $H_*(M;\bb{Z})\cong H_*(\bb{P}^1 \times 
  \dots \times \bb{P}^1;\bb{Z})$ if $M$ satisfy the hypothesis of Theorem
  \ref{th:K}. Since $H$ is 
  perfect there are exactly $\dim(H_{2k}(M))=$ $~\binom{n}{k}$ critical points of index
  $2k$. In \cite{Ha,TW},  the 
  above isomorphism  is proved by counting fixed points. We will not discuss the proof
  here.

  Denote  the points of index 2 by $p_1, \dots ,p_n$ . In the light of Theorem \ref{th:K}
  for each
  fixed point we get classes $a_1,  \dots, a_n \in H^2_{S^1}(M)$  such that
  \begin{equation}
    \begin{aligned}
      &a_j|_{p_j} = -y\\
      &a_j|_p = 0 &&\text{for all other fixed points $p$ of index 0 or 1}. 
    \end{aligned}
  \end{equation}
  
  These classes satisfy the following Proposition.

  \begin{prop}[{\cite[Prop 4.4]{TW}}]\label{p:tw}
    Let $I$ be a subset of $\{1, \dots, n\}$ with $k$ elements. There exist a unique fixed
    point $p_I$ of index $2k$ such that 
    $$a_j|_{p_I}=-y \text{  \ \ if and only if \ \ } j \in I$$
     and $a_j|_{p_I}=0$ otherwise.
  \end{prop}

  Proposition \ref{p:tw} identifies the fixed points in $M$ with subsets $I$ of
  $\mathcal{S}:=\{1, \dots, 
  n\}$. Observe that the cohomology class
  $a_I:=\prod_{i\in I}{a_i}\in H^{2k}_{S^1}(M) $ is the same as the class 
  $a_{p_I}$ mentioned in Theorem \ref{th:K}. Moreover this class is such that
  \begin{equation}\label{eq:fact}
    a_I|_{p_J}=(-1)^k y^k \text{ if and only if  }I\subseteq J
  \end{equation}
  and it is zero otherwise. 
  
  \begin{remark} 
    The class $a_0$, associated to the unique point of index zero,  takes the value $1\in
    H^0_{S^1}(pt)$ when restricted to any fixed 
    point. Therefore it is the identity element in the ring $H^*_{S^1}(M)$. Denote $ya_0$
    by $y$.   
  \end{remark}
 
  If we apply the same results to the Hamiltonian function $-H$, we obtain unique classes
  $b_J\in H^{2n-2k}_{S^1}(M)$ associated to each $p_J$ of index $2k$ such that
  $b_J|_{p_J}=(-1)^{n-k}y^{n-k}$ and is zero when restricted to all other fixed points of
  index greater 
  or equal to $2k$. These classes  also form a basis of $H^{*}_{S^1}(M)$. The next
  proposition establishes the relation with the former basis.  

  \begin{prop}\label{prop:dual}
    Let $I=\{i_1,\dots,i_k\}$ and let $I^c=\{i_{k+1}, \dots, i_n \}$ be its complement.
    Then the classes $b_I$ satisfy the following relation    
    \begin{equation}\label{eq:dual}
    b_I = \si_{n-k} + y \si_{n-k-1}+ \dots + y^{n-k},
    \end{equation}
    where $\si_i$ is the i-th symmetric function in the variables $a_{i_{k+1}}, \dots,
    a_{i_n}$.
  \end{prop}
  
  \proof
  There are two ways of seeing this. One is just by checking that when we restrict the
  right side of Equation (\ref{eq:dual}) to the fixed point $p_J$ we get by means of
  (\ref{eq:fact}) the same as $b_I|_{p_J}$. We
  can also check by hand, i.e. by restriction to all fixed points. Once again, applying
  (\ref{eq:fact}) we get  
  $$
  b_{\{i\}^c}=a_i+y,
  $$ 
  Finally we can check that 
  $$
  b_I|_{p_J}=\{\prod_{i\notin I}{(a_i + y)}\}|_{p_J} .
  $$
  for all fixed points $p_J$.

  \QED

  Consider a point $p_I$ of index $2k$ and associate  the class $a_I\in H^{2k}_{S^1}(M)$ as
  before. When we 
  restrict $a_I$ to $M$ we obtain  a class $a_I|_M\in H^{2k}(M; \bb{Z})$. By taking
  the Poincar\'e dual of $a_I|_M$, we get a homology class $\p{+}{I}\in H_{2n-2k}
  (M;\bb{Z})$. Similarly   using the class $b_I$ we get a homology class $\p{-}{I}\in
  H_{2k}(M;\bb{Z})$. Here is an immediate corollary of Proposition \ref{prop:dual}. 

  \begin{cor}\label{cor:dual}
    The class $\p{-}{I}$ is the same as the class $\p{+}{I^c}$.
  \end{cor}
  \proof This is clear because the variable $y$ is mapped to zero under reduction to
  usual cohomology. Now use that $\si_{n-k}=a_{I^c}$. 
  
  \QED

  The last part of this section establishes the relation of the $\p{\pm}{I}$ classes with
  the stable  and unstable manifolds of \S\ref{ss:morse}. This is summarized in the
  following proposition. Remember that we are working with an almost-complex structure $J$ in
  $\mathcal{J}_{\rm inv}(M)$. This result would fail without this hypothesis. 
  
  \begin{prop}\label{interplay} Let $p_I$ be a fixed point of index $2k$. Then the classes
  $\p{-}{I}$ and 
  $\p{+}{I}$ are exactly the same as the classes $[W^u(p_I)]$ and
  $[W^s(p_{I})]$ respectively.  
  \end{prop}
  \proof  
  
  Recall that $ES^1$ can be taken to be the infinite dimensional sphere $S^\infty$.
  Consider a finite dimensional  approximation
  $M^N:=M\times_{S^1}S^{2N +1}$ of $M\times_{S^1}ES^1=M\times_{S^1}S^{\infty}$ for $N\in
  \bb{N}$ big enough . These are finite dimensional smooth compact manifolds.
  Since $W^s(p_{I})$ is $S^1$-invariant, there is a natural  extension
  $W^{N,s}(p_{I}):=W^s(p_{I})\times_{S^1}S^{2N+1}$ of $W^s(p_{I})$ to $M^N$. Let $X^N$ be
  the   Poincar\'e dual of  $W^{N,s}(p_{I})$ in $M^N$.

  For all $N$, there is a natural inclusion (as fibre) $i_N:M \inject M^N$. Since the
  inclusions are natural, the restriction  ${X^N}|_M:={(i_N)}^*(X^N)\in H^*(M)$ is the same as
  the Poincar\'e  dual of  $[W^s(p_{I})]$ in $M$.

  Observe that the natural inclusions 
  $$M^N \inject M^{N+1} \inject \cdots \lim_{N}{M^N}= M\times_{S^1} ES^1$$
  induce a sequence 
  $$
   \cdots \to X^{N+2}\to X^{N+1} \to X^{N}
  $$
  given by the restrictions.  Thus, by considering the directed limit, there is an element
  $$
  X:=\lim_{N}{X^N}\in H^*(M\times_{S^1} ES^1)=H^*_{S^1}(M)
  $$
  that restricts to  $X^N$ for all $N$. Naturally, if  $i:M\inject M\times_{S^1} ES^1 $ is
  the inclusion, then  $X|_M:=i^*(X)=\PD([W^s(p_{I})])$. We claim that $X$ satisfies the
  same properties as the class $a_I$, that is, $X|_{p_I}=(-1)^k y^k$ and $X|_p=0$ for all
  other fixed points $p$ such that $\al(p)\le 2k$. Therefore, by Theorem ~\ref{th:K}
  we must have $X=a_I$. Then $\PD(X|_M)=\PD(a_I|_M)$ and the result will follow immediately. 

  Take a neighborhood  $U(p_I)$ around $p_I$ as in \S \ref{ss:morse}. Thus,
  $U(p_I)$ is isomorphic to an open neighborhood $V$ of zero in $E^+ \oplus E^-$. It is clear
  that if $U(p_I)$ is small enough, $W^s(p_{I})\cap U(p_I)$
  is diffeomorphic to  $E^+ \cap V $. Therefore, the normal bundle of $W^s(p_{I})$ can
  locally be identified with $E^-$. Finally, by carrying this localization to $X^N$ and
  considering the limit,  we have $X|_{p_I} = e(E^-)=(-1)^k y^k$, where $e(E^-)$ is the
  equivariant Euler class of $E^-$. 
  
  To finish the proof, observe that if $p$ is any other fixed point with index less than
  or equal to $2k$, there is no gradient line from 
  $p$ to $p_I$. This is because the gradient flow is Morse-Smale. Hence, by using the
  localization again we obtain that $X|_{p}=0$. This proves the proposition. 
  
  \QED  
  
  \begin{cor}\label{c:a}
    By the definition of the classes $a_I$ and $b_I$, we have 
    $$
    [W^u(p_I)]=\p{-}{I}=\PD(b_I|_M) \text{    and   } [W^s(p_I)]=\p{+}{I}=\PD(a_I|_M),
    $$
    therefore the product  $[W^u(p_I)]\cap [W^s(p_J)]$ is given by
    $$
    [W^u(p_I)]\cap [W^s(p_J)]=\PD(b_I|_M)\cap \PD(a_J|_M)=\PD(b_I a_J|_M).
    $$
  \end{cor}

  \begin{cor} By Corollary \ref{cor:dual} and Proposition \ref{interplay} above we have
    the ``duality'' relation  $[W^u(p_I)]  = [W^s(p_{I^c})]$.
  \end{cor}
  
  \begin{remark}\label{rmk:basis}
    Let $x_i:=a_i|_{M}\in H^2(M;\bb{Z})$. The theory of this section proves that 
    the elements ${x_i}$ generate the algebra $H^*(M;\bb{Z})$. Therefore a basis for
    the vector space $H^{2k}(M;\bb{Z})$ consists of the elements $x_I=x_{i_1} \dots
    x_{i_k} $ 
    for sets $I =\{ i_1< i_2 \dots < i_k \} $. Moreover, by Theorem \ref{thm:TW} the
    first 
    Chern class of $M$ is given by $c_1(M)=2(x_1 + \dots + x_n)$.
  \end{remark}

  \medskip
  
  Proposition \ref{interplay} also  provides some information about the existence of
  gradient lines. More precisely we have the next proposition.

  \begin{prop}\label{gradlines}
      Let $I=\{i_1, \dots, i_k\}\subset \mathcal{S}$. Take $i_{k+1}\notin I$ and
      consider $I'=I\cup\{i_{k+1}\}$. Let $A_I:=\sum_{i\in I}{\p{-}{i}}\in H_2(M)$. Then,   
      \begin{itemize}
      \item[a)] There is a gradient line from $p_{I'}$ to $p_I$. Moreover, the
	homology class of the sphere generated by rotating the gradient line by the
	$S^1$ action is $\p{-}{i_{k+1}}$. 
      \item[b)] There is a broken gradient line from $p_{\mathcal{S}}$ to $p_I$. The
      class 
      $A_{I^c}$ is then represented by rotating this broken line. Also,
      $\om(A_{I^c})=H_{\rm max}-H(p_I)$ and $c_1(A_{I^c})=n+m(p_{I}) $.
      \end{itemize}
  \end{prop}

  \proof
  
  To prove there is a gradient line from  $p_{I'}$ to $p_I$ we need to show that the
  intersection  $W^u(p_{I'})\cap W^s(p_{I})$ is non-empty. By definition of the
  intersection 
  product in terms of pseudocycles \cite{McS2} it is enough to prove that the intersection
  product of the 
  classes $[W^u(p_{I'})]$ and $[W^s(p_{I})]$ is non-zero. 

  Consider the equivariant cohomology classes $b_{I'}$ and $a_{I}$. By Proposition
  \ref{prop:dual} we get  
  $$b_{I'}a_{I}=a_{{I'}^c}a_I + y d$$
  where $d\in H^*_{S^1}(M)$. Since ${I'}^c=I^c \cap \{i_{k+1}\}^c=\{i_{k+1}\}^c$,
  $$
  a_{{I'}^c}a_I=a_{\{i_{k+1}\}^c}.
  $$
  Once again by Proposition  \ref{prop:dual} 
  $$a_{\{i_{k+1}\}^c}|_M=b_{i_{k+1}}|_M,$$ 
  thus
  $$
  b_{I'}a_I|_M=b_{i_{k+1}}|_M.
  $$
  Now, using Corollary \ref{c:a}  we get 
  \begin{equation}\label{eq:class}
  [W^u(p_{I'})] \cap  [W^s(p_{I})]=\PD(b_{I'}a_{I}|_M)=\PD(b_{i_{k+1}}|_M)=
  \p{-}{i_{k+1}}\ne 0. 
  \end{equation}
  Therefore, there is a gradient line, thus a whole {\em gradient sphere} $A$, just by 
  rotating the gradient line. Note that there can be more than one gradient sphere from
  $p_{I'}$ to $p_I$. We claim that all these gradient spheres must be homologous.   

  It is not hard to see from the construction of $A$ that  
  $$
  \om(A) = \int_{A}{\om} = H(p_{I'})-H(p_I).
  $$
  Therefore if $A'$ is another gradient sphere joining $p_{I'}$ and $p_I$, $\om(A)=\om(A')$.
  Also observe that if $\om'$ is any $S^1$-invariant form sufficiently close  to $\om$
  then  $\om(A)=\om'(A)$.  Now since the
  symplectic condition is an open condition we can perturb $\om$  to
  obtain a 
  new symplectic form $\om'$ close to $\om$. By averaging respect to the group action, we
  can assume the form $\om'$ to be $S^1$-invariant. This proves that the classes $A'$ and $A$
  have the same symplectic area, that is $\om'(A)=\om'(A')$, for an open set of symplectic
  forms $\om'$. Since $M$ is simply connected and there is no torsion $A$ must be
  homologous to $A'$. Finally by Equation 
  (\ref{eq:class}) this sphere must be in class $\p{-}{i_{k+1}}$. 

  To prove the second part, we can do the same process for each point in
  $I^c=\{i_{k+1}\dots i_n\}$. Then getting a sequence of gradient lines 
  
  $$
  p_{\mathcal{S}} \stackrel{\ga_1} \to p_{\mathcal{S}-\{i_{n-1}\}}  \dots  p_{I \cup
  \{i_{k+1}\}}\stackrel{\ga_{n-k}} 
  \to p_{I}. 
  $$ 
  It is clear now that the chain of gradient spheres obtained by rotating this broken
  gradient line must be in class $A_{I^c}$. Note that we could also 
  use a gluing 
  argument as in \cite{Sch2} to prove that there is an honest gradient line from
  $p_{\mathcal{S}}$ to $p_I$.   Thus $\om(A_{I^c})=H_{\rm max}-H(p_I)$ and
  $c_1(A_{I^c})=m(p_I)-m(p_{\mathcal{S}})=n+m(p_I)$.

  \QED


  \section{Quantum Cohomology and the Seidel Automorphism}

  \subsection{Small Quantum Cohomology}\label{ss:qh}
  In the literature, there are several definitions of quantum cohomology. In this section 
  we make precise the definition of the quantum cohomology we are using,  assuming the
  definition of   genus zero Gromov-Witten invariants. We will follow entirely the
  approach of \cite[Chapter 11]{McS2}. 

  Let $\La_\om$ be the usual {\em Novikov ring} of $(M,\om)$. We recall that $\La_\om$ is the
  completion 
  of the group ring of $H_2(M):=H_2(M;\bb{Z})/{\rm Torsion}$. It consists of all (possibly
  infinite) formal sums of the form 
  $$
  \la=\sum_{A\in H_2(M)}{\la_A e^A}
  $$
  where $\la_A\in \bb{R}$ and the sum satisfies the finiteness condition
  $$
  \#\{ A\in H_2(M) | \la_A \ne 0, \om(A) \le c \} < \infty
  $$
  for every real number $c$. By definition, $\deg(e^A)=2c_1(A)$, where $c_1$ is the first
  Chern class of $M$.

  The {\bf (small) quantum cohomology } of 
  $M$ with coefficients in $\La_\om$ is  defined by  
  $$ 
  QH^*(M):=H^*(M)\otimes_{\bb{Z}}\La_\om.
  $$
  As before $H^*(M)$ denotes the ring $H^*(M;\bb{Z})$ modulo torsion.
  We now proceed to define  the {\bf quantum product} on $QH^*(M)$. We want the quantum
  product to be a linear homomorphism of $\La_\om$-modules 
  $$
  QH^*(M)\otimes_{\La_\om} QH^*(M) \to QH^*(M) : (a,b)\mapsto a*b.
  $$
  Since  $QH^*(M)$ is generated by the elements of $H^*(M)$ as a $\La_\om$-module, it is
  enough 
  to describe the multiplication for elements  in $H^*(M)$.    
  Let $e_0, e_1, \dots , e_n$ be a basis for $H^*(M)$ (as a $\bb{Z}$-module). Assume each
  element is homogeneous and $e_0=1$, the identity for the usual product. Define the
  integer matrix 
  $$
  g_{ij}:=\int_{M}{e_i\smile e_j}.
  $$
  Here $e_i\smile e_j$ is the usual cup product in cohomology. Let $g^{ij}$ be the inverse
  matrix.  The quantum product of  $a, b\in H^*(M)$, is defined by
  \begin{equation}\label{eq:qprod}
  a*b:=\sum_{B\in H_2(M)}\sum_{k,j} {\GW{B}{3}(a,b,e_k)g^{kj}e_j \e{B}}.
  \end{equation}
  The coefficients $\GW{B}{3}$ are the usual Gromov-Witten invariants of
  $J$-ho\-lomorphic curves in class $B$. The terms in the sum are nonzero only if $\deg
  (e_k) + \deg (e_j) = \dim M$ and  $\deg (a)+\deg (b) +\deg (e_k) = \dim M +
  2c_1(B)$. Thus, it is enough to consider classes $B$ such that 
  $$
  \deg (a) +\deg (b) - \dim M\le 2c_1(B) \le \deg(a) + \deg (b).
  $$

  In the problem at hand, a basis for $H^*(M)$ is given by the elements $x_I$ as in
  \ref{rmk:basis}. Then the integrals  
  $$
  g_{IJ}=\int_{M}{x_I \smile x_J}
  $$
  all  vanish unless the sets $I$ and $J$ are complementary. This is because  if $I$,
  $J\subset\{ 1, \dots, n\}$,  $x_I\smile x_J=\xp $ if and only if $I^c=J$. Here $\xp$ is the
  positive generator of $H^{2n}(M;\bb{Z})$. 
  
  We claim that to compute the quantum
  product, we only need to consider in Equation (\ref{eq:qprod}) classes $B$ such that
  $c_1(B)\ge 0$. More precisely, we have the proposition.
  
  \begin{prop}\label{jump} Assume $(M,\om)$ is a symplectic manifold with a
    semi-free $S^1$-action with only isolated fixed points.
    Let $B\in H_2(M)$, and  let $a,b,c\in H^*(M)$. 
    If $c_1(B)<0$, then the Gromov-Witten invariant $\GW{B}{3}(a,b,c)$ is zero. Moreover,
    if $c_1(B)=0$ and some $\GW{B}{3}\ne 0$, then $B=0$. Therefore,  
    the expression for the quantum product (\ref{eq:qprod}) can  be written as  
    $$
    a*b=a\smile b + \sum_{B\in H_2(M), c_1(B)> 0} {a_{B} \e{B}}.
    $$
    where the classes $a_{B}$ have degree $\deg( a_{B} )=\deg a + \deg b - 2c_1(B)$.
  \end{prop}
  
  \begin{rmk}\label{rmk:jump}
    Note that since $c_1(B)$ is even, the classes $a_{B}$  appear in the sum above by
    ``jumps'' of four in the degree.  
  \end{rmk}
  
  The rest of this section is dedicated to the  proof of Proposition ~\ref{jump}.
  
  \medskip
  
  To compute the Gromov-Witten invariants $\GW{B}{3}(a,b,c)$ one usually constructs a
  regularization 
  (virtual cycle) $\vmc(M,J,B)$ of the moduli space $\moduli(M,J,B)$. Then one computes
  the  
  intersection number of the evaluation map 
  $$ev:\vmc(M,J,B)\to M^3$$ 
  with a cycle $\ag_1\times \ag_2 \times \ag_3$  representing the class
  ${\rm PD} (a) \times {\rm PD} (b) \times {\rm PD} (c) $. This procedure can be modified
  in the  
  following way. First, let $\ag:Z \to M^3$ be a pseudocycle that represents the product 
  ${\rm PD} (a) \times {\rm PD} (b) \times {\rm PD} (c) $, then  define the {\em cut-down}
  moduli space  by 
  $$\moduli(M,J,B;Z):=ev^{-1}(\overline{\ag(Z)}).$$
  Here $ev: \moduli(M,J,B) \to M^3$ is 
  the evaluation map and $\overline{\ag(Z)}$ is the closure in $M$ of the
  pseudocycle $Z$ \cite{McS2}. Finally, construct a regularization of the cut-down moduli
  space. McDuff and Tolman use this approach to calculate the Gromov-Witten 
  invariants. The next two results are proved in \cite{McT}. They show
  exactly how to compute the invariants $\GW{B}{3}$ using this procedure. Remember that an
  $S^1$ action on $M$ can be extended to 
  an action on $J$-holomorphic curves just by post-composition. Also, a  pseudocycle
  $\ag:Z \to M$ is said to be $S^1$-invariant, if $\ag(Z)$ is.  
  
  \begin{prop}\label{prop:dusa2} 
    Let $(M,\om)$ be a symplectic manifold. Then, the Gromov-Witten invariant
    $\GW{B}{3}(a,b,c)$ is a sum of contributions, one 
    from each connected component of the moduli space 
    $\moduli(M,J,B;Z)$.    
    
    Assume now that $M$ is equipped with an $S^1$ action $\{\la_t\}$, and that
    $\ag:Z \to M^3$ and $J$ are $S^1$-invariant. Then, a connected component of
    $\moduli(M,J,B;Z)$ makes no contribution to $\GW{B}{3}(a,b,c)$ unless it contains an
    $S^1$-invariant element. 
      
  \end{prop}

  The following  lemma describes what  the invariant elements in the moduli space
  $\modu{k}(M,J,B)$ are. We include a proof so that Corollary \ref{c:inv} is a more natural
  result. 
  \begin{lemma}\label{l:inv} 
    Let $(M,\om)$ be a symplectic manifold with a semi-free
    $S^1$-action.  Let $[u]$ be a class in the moduli space $\modu{k}(M,J,B)$ represented
    by a $J$-holomorphic sphere $u:\bb{P}^1 \to M$. Assume $[u]$ is fixed by the action
    $\la = \{ 
    \la_\tg\}$. Then, there are at most two marked points, i.e. $k\le2$ and $u$ can be
    parametrized as   
    $$
    u:\bb{R} \times S^1 \to M, \ \  u(s,t)=\la_{pt}\ga(s).
    $$
    Here $\ga:\bb{R}\to M$ is a path joining two fixed points $x,y\in M$ so that the
    marked points are in $u^{-1}\{ x,y\}$, and $\ga$ satisfies the
    gradient flow 
    equation 
    \begin{equation}\label{eq:grad}
      \ga'(s)=p\  \grad(H) \text{ \ for some \ }p\ne 0.
    \end{equation} 
    Moreover, if we fix $\ga$, the parametrization is
    unique provided
    $$
    x=\lim_{s\to -\infty}\ga(s) \text{  and  } y=\lim_{s\to \infty}\ga(s).
    $$
  \end{lemma}
  
  \proof
  Let $u:\bb{P}^1 \to M $ be a non constant and not multiply covered $J$-holomorphic sphere
  in $M$. For each $\tg\in S^1$ the map $\la_\tg \circ u$ must be a reparametrization of
  $u$. This 
  is because the equivalence class $[u]$ is fixed under the action. Thus,
  there is a $\phi_\tg\in  \rm{PSL}(2,\bb{C})$ such that  $ \la_\tg \circ u = u \circ
  \phi_\tg$. Since the map 
  $u$ is not multiply covered $\phi_\tg$ is unique. Then, it is easy to see that the
  assignment $S^1 \to \rm{PSL}(2,\bb{C}) : \ \tg \mapsto  \phi_\tg $ is a
  homomorphism. Since the only circle subgroups of  
  $\rm{PSL}(2,\bb{C})$  are rotations about an axis, we can choose coordinates on $\bb{P}^1$
  so that the rotation axis is the line joining the unique fixed points $[0:1]$ and
  $[1:0]$. Assume that ${\rm Im}(u)\cap{M^{S^1}}=\{x,y\}$. Identify $\bb{P}^1 / \{
  [0:1],[1:0] \}$ with the cylinder $\bb{R} \times S^1$ 
  with complex structure $j_0$ defined by $j_0(\partial_s)=\partial_t$, $(s,t)\in \bb{R}
  \times S^1$.  If $k=2$ we 
  identify the marked points $[0:1],[1:0]$ with the ends of the cylinder, so that
  $u([0:1])=x$ and $u([1:0])=y$. In general the image of the marked points must be fixed
  by the action. Therefore the marked points can be identify with a
  subset of $\{[0:1],[1:0]\}$. If  $(s,t)\in \bb{R} \times S^1$ are the standard
  coordinates, then 
  $$
  \phi_\tg(s,t) = (s, t +q \tg ), \text{  and } (\la_\tg \circ u)(s,t) = u (s, t +q \tg )
  \text{  where } q=\pm 1.
  $$
  
  Define $\ga(s):=u(s,0)$. Then we get $u(s,t)=\la_t\ga(s)$. Since $u$ is
  $J$-homomorphic and $J$ is invariant
  $$(\la_t)_* (\ga'(s) + JX(\ga(s)))=\partial_s u + J\partial_t u=0.$$
  
  With respect to the metric $g_J$, the gradient flow of $H$ is given by $\grad H=-JX$, thus
  $\ga'(s)=\grad(H)(\ga(s))$. Now use that any sphere is a $|p|$-fold cover of a simple
  one. We absorb any negative sign into $p$ rather than $q$.  
  \QED
  
  Note that our original goal was to understand the invariant  stable maps in
  $\moduli(M,J,B;Z)$. By the 
  previous lemma, the non-constant components of the stable maps may carry at most two
  special 
  points. Then the $S^1$-invariant elements in $\moduli(M,J,B;Z)$ may have a {\em ghost}
  component that carries the third marked point.

  We have as an  immediate consequence the following corollary. 
  
  \begin{cor}\label{c:inv} Assume the same hypothesis as in Lemma \ref{l:inv}.
    Let  $u$ be an $S^1$-invariant sphere, and let $A\in H_2(M)$ be its homology class in
    $M$. Then its 
    first Chern class is given  by $c_1(A)=p(m(x)-m(y))$, and  is always
    positive. Here $m(x)$ is the sum of the weights at $x$.
  \end{cor}
  
  \proof
  Recall that $m(x)=n-\al(x)$ where $\al(x)$ is the Morse index of  $x$. Now, 
  if $m(x)<m(y)$, the path $\ga$ from $x$ to $y$ must satisfy  Equation (\ref{eq:grad})
  with $p<0$. This is because there are no generic solutions to this 
  equation  otherwise. Then $c_1(A)=p(m(x)-m(y))$. If  $m(x)>m(y)$, now $p$ must be
  positive and the result follows.
  
  \QED
  
  \begin{rmk}\label{r:positive}
    Let $u$ be an $S^1$-invariant holomorphic sphere, let   $A\in H_2(M)$ be its homology
    class. Lemma \ref{l:inv} and Corollary \ref{c:inv} imply that if  $c_1(A)=0$ then $A$
    must be zero. This is because if $A$  joins two fixed points
    $x,y\in M$, they  must have the same index, which is not possible because the flow is
    assumed to be Morse-Smale. 
  \end{rmk}
  
  \proof[Proof of Proposition \ref{jump}]
  By Proposition \ref{prop:dusa2} a component of $\moduli(M,J,B;Z)$  contributes to 
  $\GW{B}{3}(a,b,c)$ only if the moduli space  has a $S^1$-invariant stable map ${\bf
    u}$. We can assume that there is at least one non-trivial component  
  $u_i$ of the stable map ${\bf u}$. Since ${\bf u}$ is invariant, so is
  $u_i$. Therefore,  Corollary \ref{c:inv} implies that $c_1(B_i)>0$ if $B_i \in
  H_2(M)$ is the class of $u_i$.  
  Then $c_1(B)>0$ and  the first claim follows. Note that the second part is
  a direct consequence Lemma \ref{l:inv}, because any  $S^1$-invariant $J$-holomorphic map with zero
  Chern class must be constant.   
  
  Finally , the product $a*b$ can  be written as  
  $$
  a\smile b + \sum_{c_1(B)>0}\sum_{I} {\GW{B}{3}(a,b,\x{I})\xc{I} \e{B}}.
  $$
  Now take 
  $$ 
  a_B:=\sum_{I} {\GW{B}{3}(a,b,\x{I})\xc{I}}.
  $$
  This proves the proposition. Note that we have $\deg( a_{B} )=\deg a + \deg b - 2c_1(B)$.
  
  \QED

  \subsection{Almost Fano Manifolds}\label{afano}
  Assume the hypothesis of Proposition \ref{jump}. The relevant spheres (the ones that
  count for the  GW invariants) all have  
  positive first Chern class.  Moreover, let $B\in H_2(M)$ be as in Proposition
  \ref{jump}, then $c_1(B)>0$. Since
  $B$ is invariant, using Proposition  \ref{gradlines} and Lemma
  \ref{l:inv}, $B$ can be written as a combination 
  $$
  B=\sum_{i}{d_i \p{-}{i}}
  $$
  where the coefficients $d_i$ are non-negative integers.  Therefore, if we define
  $A_i:=\p{-}{i}$ and $q_i:=e^{A_i}$,  we may now consider the  polynomial ring 
  $$\La=\bb{Q}[q_1, \dots, q_n]$$ 
  as coefficients for the quantum cohomology. Then, if $B$
  is as before, 
  $$e^B={q_1}^{d_1}\dots {q_n}^{d_n}.$$ 
  This will be really useful  in \S
  \ref{ss:proof}. For the rest of this paper, we will  assume $\La$ to be the quantum
  coefficient ring.
  
  We finish this section with a discussion about the behavior of $J$-holomorphic curves
  in $M$. In the literature an almost complex manifold $(N, J)$ is said to be {\bf Fano }
  if the first Chern class $c_1(TN,J)$ takes positive values on the {\bf effective cone}
  $\eff(N,J)$, namely
  $$\eff(N,J):=\{ A \in H_2(N) | \exists \text{ a $J$-holomorphic curve in class } A \}.$$
  
  In symplectic geometry sometimes is useful to consider the definition
  
  $$\eff(N,\om)=\{ A \in H_2(N)|  A_1,\dots , A_n \in H_2(N): A=\sum_i{A_i}, \GW{A_i}{3}
  \ne 0\}$$ 
  for the effective cone on a symplectic manifold $(N,\om,J)$ with a compatible almost
  complex structure $J$. Its clear that $\eff(N,\om) \subset \eff(N,J)$. Then, we can say
  that $(N,\om,J)$ is {\bf almost Fano} 
  if the first Chern class $c_1(TN,J)$ takes positive values on the effective cone
  $\eff(N,\om)$. We have the following corollary.

  \begin{cor}
   Let $(M,\om)$ be a symplectic manifold with a semi-free $S^1$-action with isolated
   fixed points. Then $(M,\om,J)$ is almost Fano.
  \end{cor}
  

  \medskip

  \subsection{The Seidel Automorphism}\label{ss:seidel}
  
  In this paragraph we introduce the theory behind the definition of the Seidel element.  
  The results concerning the present problem are discussed next. We will follow closely the 
  book \cite{McS2}. The proofs of the results exposed in this section are mostly
  contained in  Chapters 8,9 and 11. 

   Let $M$ be as in \S \ref{intro}. Since the action is Hamiltonian, it is possible to
  associate to $M$ the locally trivial bundle $\T{M}_\la$ over $\bb{P}^1$ with fibre $M$
  defined by the {\em clutching function} (action) 
  $\la:S^1 \to {\rm Ham}(M,\om) $: 
  $$
  \T{M}_\la:=S^3  \times_{S^1} M
  $$
  
  We denote the fibres at $[1:0]$ and $[0:1]$ by $M_0$ and $M_\infty$ respectively. 
  Note that the isomorphism type of  $\T{M}_\la$ only depends on the homotopy class of
  $\la$. 

  Since $\la$ is Hamiltonian, we can construct a symplectic form $\Og$ on
  $\T{M}_\la$. In fact the bundle $\pi: \T{M}_\la \to \bb{P}^1$ is a {\em Hamiltonian
  fibration} with fibre 
  $M$, thus admitting sections (\cite[Chapter 8]{McS2}). 

  In the case when the manifold has an $S^1$-action,  we choose an $\Og$-compatible
  almost complex structure $\T{J}$ on $\T{M}$, such that $\T{J}$ is the product $J_0
  \times J$ 
  under trivializations. We can define for each fixed point $x\in M^{S^1}$ a
  pseudoholomorphic  section $\si_x:=\{[z_0:z_1;x]:[z_0:z_1]\in \bb{P}^1\}$.

  Take $\T{A}\in H_2(\T{M}_\la, \Og)$ a section class, that is $\pi_*
  (\T{A})=[\bb{P}^1]$. Let $a_1, a_2 \in H^*(M)$. Given two fixed marked points $w_1,
  w_2\in \bb{P}^1$ we may think of the Poincar\'e dual to the class $a_i$ as
  represented by a cycle $Z_i$ in the fibre  $M_i \inject \T{M}_\la$ over $w_i$. With this
  information it is possible to construct the Gromov-Witten invariant
  $\GW[\T{M}_\la, {\rm \bf w}]{\T{A}}{2}(a_1, a_2)$.  This invariant counts
  the number of $J$-holomorphic sections of $\T{M}_\la$ in class $\T{A}$ that pass through
  the cycles $Z_i$.

  \begin{definition}
    Let $(M, \om)$ be as before. Let $\si:\bb{P}^1\to \T{M}_\la$ be a section. The {\bf Seidel
      automorphism} 
    $$
    \sei{\si}:QH^*(M;\La)\to QH^*(M;\La)
    $$
    is defined by 
    \begin{equation}
     \sei{\si}(a)= \sum_{A\in H_2(M)} \sum_{k,j} {\GW[\T{M}_\la, {\rm \bf w
     }]{[\si]+i_*A}{2}(a,e_k)g^{kj}e_j}\e{A}.     
    \end{equation}
    where $i:M\to \T{M}_\la $ is an embedding (as fibre). 
  \end{definition}
  In this definition we are considering  a  basis $\{e_i\}$ for $H^*(M)$ as in Equation
  ~(\ref{eq:qprod}). It is important to remark that the Seidel automorphism as defined
  above does not preserve 
  degree. The shift on the  degree depends  on the section class $\si$ that we use as
  reference. 

  If $\1 \in QH^*(M)$ denotes the identity in the quantum cohomology ring, the class
  $\sei{\si}(\1)\in QH^*(M)$ is  called the {\bf Seidel Element} of the action respect to
  the section $\si$. We will use the same notation for the Seidel automorphism and the
  Seidel element. Thus,  the Seidel automorphism
  is now given just by quantum multiplication by  the element $\sei{\si}$  \cite{McS2}. That is, 
  $$
  \sei{\si}(a) = \sei{\si}*a.
  $$
  Note that the Seidel automorphism shifts degree by  $\deg(\sei{\si})$.


  \subsection{Seidel Automorphism and Isolated Fixed Points}\label{ss:seidel2}

  Consider now the present problem. That is, assume that the action is semi-free and it has
  isolated fixed points.  Let $\simax$ be the section defined by the fixed
  point $p_\mathcal{S}$. In this particular case the automorphism  $\sei{\simax}$ increases the
  degree by $2n$. Let $p_I\in M$ be a fixed point. Recall that  we can associate to $p_I$
  classes in homology $\p{-}{I}$ and $\p{+}{I}$, and if we consider all the fixed points, then
  the classes $\p{+}{I}$ form a basis for $H^*(M)$.  

  The next theorem, due to McDuff and
  Tolman \cite{McT}, gives the first step towards a description of the Seidel
  automorphism. Although they have proved this result in great generality (the fixed
  points are  allowed to be in submanifolds rather than being isolated) and they use
  quantum homology  rather than cohomology, it is not hard to adapt their result  to our
  present notation.     

  \begin{thm}[McDuff-Tolman] Let $(M,\om)$ be a symplectic manifold with a
    semi-free circle action with isolated fixed points. Assume its associated Hamiltonian 
    function $H$ is such that $\int_{M}H\om^n=0$. Let $A_I\in H_2(M)$ be as
    considered in \ref{gradlines}. Then, the Seidel automorphism can be expressed as    
    $$\sei{\simax} (\PD(\p{-}{I}))=\PD(\p{+}{I}) \e{A_{I^c}} + \\ 
    \sum_{\om(B)>0} {a_{B} \e{A_{I^c} + B}}.$$
    where $a_B\in H^*(M)$. If  $a_B \ne 0$ then  $\deg{x_I}-\deg{a_B}=2c_1(B)$. 
    Moreover, if we write the sum above in terms of the basis $\{\PD(\p{+}{J})\}$ we get
    $$\sei{\simax} (\PD(\p{-}{I}))=\PD(\p{+}{I})\e{A_{I^c}} + 
    \sum_{\om(B)> 0} \sum_{J\in \mathcal{S}} {C_{B,J}  \ \PD(\p{+}{J})\e{A_{I^c} + B}}. 
    $$ 
  \end{thm}

  We know by Corollary \ref{cor:dual} that $\p{-}{I}=\p{+}{I^c}$. By definition
  $\PD(\p{+}{J})=x_J$, therefore  we have the following straightforward corollary.

  \begin{cor}\label{contri} Let $(M,\om)$ be a symplectic manifold with a
  semi-free circle action with isolated fixed points. Assume its associated Hamiltonian
  function $H$ is such that $\int_{M}H\om^n=0$. Let $\{x_I\}$ be the basis for the
  cohomology ring as considered in Remark \ref{rmk:basis}, and let $A_I\in H_2(M)$ as
  considered in \ref{gradlines}. The Seidel automorphism can be expressed as
  \begin{equation}\label{eq:seia}
  \sei{\simax}(\xc{I}) = \x{I} \e{A_{I^c}} + 
  \sum_{\om(B)> 0} \sum_{J\in \mathcal{S}} {C_{B,J}  \ \x{J}\e{A_{I^c} + B}}. 
  \end{equation}
  The rational coefficients $C_{\T{B},J}$ can be nonzero only if $|I|-|J|=c_1(B)$ and
  the moduli space $\mo{\overline{W^u(p_I)}}{\overline{W^u(p_J)}}$
  has  an $S^1$-invariant element. $\si_I$ denotes the section defined by the fixed
  point $p_I$. 
  \end{cor}   
  
  Thus, the key to understand the Seidel automorphism is first
  to know what  the $S^1$-invariant elements in moduli spaces 
  $$\mo{Z}{Z'}$$ are. Here $Z$ and $Z'$ are closed $S^1$-invariant  cycles in $M$. These 
  elements are called {\em invariant chains in section class $\si_z + A$ from  
  $x\in Z$ to  $y\in Z'$ with root $z$} \cite{McT}. We will explain what is the meaning of
  this.   
  
  Given $x,y,z\in M^{S^1}$ an invariant {\em principal} chain in section class $\si_z + A$
  from   
  $x\in Z$ to  $y\in Z'$ with root $z$ is a sequence of fixed points
  $x=x_1, \dots, x_k=y$ joined by $\T{J}$-holomorphic spheres with the following properties:
  \begin{itemize}
  \item[a)] There is $1\le i_0\le k$ such that $x_{i_0}=x_{i_0 +1}=z$, and they are joined
    by the section $\si_z$.
  \item[b)] For each $1\le i< k$ where $i\ne i_0$, the points $x_i, x_{i+1}$ are joined
    by an invariant sphere (in $M$) in class $A_i$.
  \item[c)] $\sum_{i\ne i_0}A_i = A$.
  \end{itemize}
  An {\bf invariant chain} in section class $\si_z +A$ from  
  $x\in Z$ to  $y\in Z'$ with root $z$ is a chain as above with additional ghost
  components at each of which  a tree of invariant spheres is attached. In this case, $A$
  is the sum of classes represented by the principal spheres and the bubbles.

  Also, we can decompose $A=A' + A''$, where $A'$ is the sum of spheres embedded in 
  the fibre $M_0$ and $A''$ the ones in $M_\infty$. An immediate lemma is the following

  \begin{lemma}\label{ineq} Assume the hypothesis of Corollary \ref{contri}, and suppose
  $\si_z + A$ is an invariant chain in the moduli space  
  $$
  \mo{\overline{W^u(p_I)}}{\overline{W^u(p_J)}}.
  $$
  Let $A=A'+A''$ be the decomposition of $A$ as described above. Then, the first Chern
  classes $c_1(A'), c_1(A'')$ can be estimated by  
  $$c_1(A')\ge |m(x) - m(z)| \text{ and } c_1(A'')\ge |m(y) - m(z)|.$$
  Therefore 
  \begin{equation}\label{eq:estimate}
    \begin{aligned}
      c_1(A)&\ge |m(x) - m(z)| + |m(y) - m(z)|, \\
      c_1(B)&\ge \max\{c_1(A'), c_1(A'')\}.
    \end{aligned}
  \end{equation}
  Moreover if the coefficient
  $C_{B,J}\ne 0$, then $c_1(B)> 0$. Finally, observe that $c_1(A)=0$ if
  and only if  $A=0$. 
  \end{lemma} 

  \proof If $A_i$ is an invariant sphere joining $x_i$ to $x_{i+1}$, Lemma \ref{l:inv}
  shows 
  that $c_1(A_i)\ge |m(x_i) - m(x_{i+1})|$ . Then
  $c_1(A')\ge\sum_{i=0}^{k}{|m(x_i)-m(x_{i+1})|}\ge |m(x)-m(z)|$. The other part is
  analogous. Now, write $\si_I + B= A + \si_z$, since $x\in W^u (p_I)$, $m(x)>m(p_I)$, then
  $c_1(B)\ge c_1(A'')$. Similarly  $c_1(B)\ge c_1(A')$. For the last statement, note that
  if $C_{B,J}\ne 0$ then $A\ne 0$. Then $A'\ne 0 $ or $A''\ne 0 $. In any case
  $c_1(B)> 0$. For the last claim, note that if $A_i$ is an invariant sphere with
  $c_1(A_i)=0$, Remark \ref{r:positive} implies that $A_i$ must vanish.
  
  \QED

  With Lemma \ref{ineq} we can simplify the expression (\ref{eq:seia}) to get the
  following corollary.

  \begin{cor}\label{contri-b}
    Assume the same hypothesis of Corollary \ref{contri}. Then the Seidel element is given
    by 
    \begin{equation}\label{eq:seib}
      \sei{\simax}(\xc{I}) = \x{I} \e{A_{I^c}} + 
      \sum_{\om(B)> 0, c_1(A)>0} \sum_{J\in \mathcal{S}} {C_{B,J}  \ \x{J}\e{A_{I^c} + B}}.  
    \end{equation}
    Again  $C_{B,J}=0$ unless $|I|-|J|=c_1(B)$
    and the moduli space $\mo{\overline{W^u(p_I)}}{\overline{W^u(p_J)}}$
    has  an  $S^1$-invariant element. 
  \end{cor}
  
  Note that the only difference to Equation (\ref{eq:seia}) is that we are considering
  only classes $B$ with positive Chern number.
  
  If there are any  higher order terms, that is, terms that correspond to
  positive first Chern classes $c_1(B)>0$, they contribute  to the sum (\ref{eq:seib}) as
  an element of degree $2( |J| + c_1(A_{I^c} + B))$. Heuristically an  invariant chain $A +
  \si_z$  makes a contribution only if $c_1(A)$ is big enough so that the inequalities
  (\ref{eq:estimate}) are satisfied.  We will see
  in our next  result that with our present hypotheses there are no such
  contributions. Thus there are not higher order terms. This result fails if for instance
  we allow the action to have fixed points along submanifolds, as we will see in the example
  described in  \S \ref{ex:blowup}.  
  Observe that we can normalize our Hamiltonian function $H$ (by adding a constant) so
  that $\int_{M}{H{\om}^n=0}$ without altering any of our previous results.    

  \begin{thm}\label{seidelact}Let $(M,\om)$ be a symplectic manifold with a
  semi-free circle action with isolated fixed points. Assume its associated Hamiltonian
  function $H$ is such that $\int_{M}H\om^n=0$. Then, the Seidel automorphism $\sei{\simax}$
  acts on the basis $\{ \x{I} \}$ by  
    \begin{equation}
      \sei{\simax} (\x{I})=\xc{I} \e{A_I}
    \end{equation}    
  \end{thm}
  
  \proof 
  Consider $I^c$ instead of $I$. By Corollary \ref{contri-b} the Seidel automorphism can
  be computed
  $$ 
    \sei{\simax} (\xc{I})  = \x{I} \e{A_{I^c}}\\ 
    + \sum_{c_1(B)>0, J\in \mathcal{S}} {C_{B,J} \x{J} \e{A_{I^c} +B} }
  $$
  
  As in Proposition \ref{jump}, the Chern number  $c_1(B)$ is a multiple of
  two. Thus the terms in the sum appear with  ``jumps'' of four in the degree. By
  Corollary 
  \ref{contri-b}, $C_{B, J}$ is nonzero only if there is a $S^1$-invariant 
  element in the moduli space $\mo{\overline{W^u(p_I)}}{\overline{W^u(p_J)}}$. We want to 
  see that the coefficients  $C_{B,J}$ are all zero.

  By contradiction assume there is an invariant chain $\si_z + A$  in this moduli space. 
  Therefore $A$ goes from  a fixed point $x\in \overline{W^u(p_I)}$  to a fixed point $y
  \in \overline{W^u(p_J)}$. This chain satisfies 
  \begin{equation}\label{eq:chain}
    \si_z + A = \si_I + B.
  \end{equation}
  
  Since the gradient flow is Morse-Smale and there is a gradient line from $p_I$ to $x$,
  $m(x)\ge m(p_I)=n-2|I|$.  Analogously $m(y) \ge m(p_J)=n-2|J|$. Since   $    c_1(B)=|I|
  -  |J|>0$ and we know $c_1(A)+m(z)=m(p_I)+c_1(B)$ from Equation (\ref{eq:chain}), we get 
  \begin{equation}\label{eq:A}
    c_1(A)=2|K|-|I|-|J|,
  \end{equation}
  where $K \subset \mathcal{S}$ is such that $p_K=z$.

  Finally, from Lemma \ref{ineq} we have 
  
  $$
  \begin{aligned}
    c_1(A) &\ge  |m(x) - m(z)| + |m(y) - m(z)|  \\
    &\ge -2m(z) +  m(y)+m(x)\\
    &\ge 4|K| - 2|I| -2|J|.
  \end{aligned}
  $$ 
  
  Therefore, by Equation (\ref{eq:A}) 
  $$
  2|K|-|I|-|J|=c_1(A)\ge 2(2|K|-|I|-|J|).
  $$
  This is only possible if $c_1(A)=0$, i.e
  $2|K|-|J|=|I|$. By Lemma \ref{ineq} $A$ must be zero. Thus
  $x=y=z$. Therefore  $B= \si_z - \si_I$. Hence
  $c_1(B)=m(z)-m(p_I)=2(|I|-|K|)$. Since $c_1(A)=0$, Equation (\ref{eq:A}) implies
  $|I|-|K|=|K|-|J|$. Thus  $0<c_1(B)=2(|K|-|J|)$. By hypothesis  $p_K=z=y \in
  \overline{W^u(p_J)}$. Then  we have $|K|\le|J|$. Thus $c_1(B)\le 0$, which is a
  contradiction. This proves the theorem.    
  
  \QED
  
  \begin{cor}\label{cor:prod}
    The Seidel element $\sei{\simax}$ is given by 
    $$
    \sei{\simax}=\xp.
    $$ and the quantum product of $\xp$ with the element $\x{I}$ is given by
    \begin{equation}\label{eq:action2}
      \xp*\x{I}=\xc{I}\e{A_I}. 
    \end{equation}   
  \end{cor}
  \proof
  The first part is obvious since $\sei{\simax} =  \sei{\simax}*\1=\sei{\simax} *
  x_0=\xp\e{0}$.  For the second part, observe that  
  $$
  \xc{I} \e{A_I} =\sei{\simax} *\x{I}=\xp * x_I.
  $$
  \QED
  
  The next paragraph is dedicated to discuss an  example where the symplectic manifold has a
  semi-free circle action but  the Seidel automorphism has higher order terms when
  evaluated on a particular class. In this example the fixed points are along 
  submanifolds. This illustrates that we cannot have a result similar to Theorem
  ~\ref{seidelact} if we weaken one of our hypothesis. 
  \subsection{Example}\label{ex:blowup} \cite[Example 5.1]{McT}
  Let  $M = \T{\bb{P}^2}$ be the one point blow up of $\bb{P}^2$ with the symplectic
  form $\om_\mu$ so that on the exceptional divisor $E$,  $0<\om_\mu (E)=\mu<1$ and
  if $L=[\bb{P}^1]$ is the standard line, we have $\om_\mu(L)=1$. We can identify $M$ with
  the space
  $$
  \{ (z_1,z_2)\in \bb{C}^2 | \ \mu\le |z_1|^2 + |z_2|^2 \le 1 \}
  $$
  where the boundaries are collapsed along the Hopf fibres. One of the collapsed
  boundaries is identified with the exceptional divisor. The other with $L$.

  A basis for $H_*(M)$ is given by the class of a point $pt$, the exceptional divisor
  $E$, the fibre class $F=L-E$ and the fundamental class $[M]$. Note that the intersection
  products are given 
  by $E \cdot E= -1$, $E\cdot F = 1$, $F\cdot F=0$. Denote by $b$ and $f$ the Poincar\'e
  duals of $E,F$ respectively. Then $b\cdot b=-1 $ and $f\cdot f=0$. It is not hard to
  see that the positive generator of $H^4(M)$ is $b\smile f=\PD(pt)$. Let us denote this
  class by just $bf$, so that a basis for the cohomology ring is $\{\1, b, f, bf \} $. Observe that
  $M$ with the usual complex structure is Fano .
  
  The non-vanishing Gromov-Witten invariants are given by
  $$
  \begin{aligned}
  \GW{L}{3}(bf,bf, f)&=\GW{F}{3}(bf, b,b)=1;\\
  \GW{E}{3}(c_1,c_2,c_3)&=\pm 1 \text{   where  } c_i = b \text{ or  }f.
  \end{aligned}
  $$
  Let us consider the usual Novikov ring $\La_\om$ as the quantum coefficients. Then the
  quantum products are give by:
  $$
  \begin{array}{lclcclcl}
    bf*bf&=& (b+f)\e{L}  &&&  bf*f&=&\1\e{L}\\
    bf*b &=& f\e{F}      &&&  b*b&=&-bf + b\e{E} + \1\e{F}\\
    b*f  &=& bf -b\e{E}  &&&  f*f&=&b\e{E}.  
  \end{array}
  $$
  
  In \cite{McT} it is proved that the circle action on $M$ given by:
  $$
  \ag : (z_1,z_2)\mapsto (e^{-2\pi it}z_1,e^{-2\pi it}z_2 ), \ \ \ \ \ \text{for }0\le t \le 1.
  $$
  is Hamiltonian. The maximum set of this action is exactly the points lying on the
  exceptional divisor $E$ and the minimum set is the line $L$. After taking an appropriate
  reference section $\si$, the Seidel element  $\Psi(\ag, \si)$ is given by 
  $$
  \Psi(\ag, \si)=b.
  $$
  Thus, evaluating the Seidel map on the class $f$ we have 
  $$
  \Psi(\ag, \si)(f)=\Psi(\ag, \si)*f=b*f= bf -b\e{E}.
  $$
  Therefore the Seidel automorphism does have higher order terms when evaluated on the
  class $f$.

  \medskip
  \section{Proof of main result}\label{ss:proof}  
  Now we are ready for proving the main theorem. Recall that the quantum coefficient ring
  is $\La=\bb{Q}[q_1, \dots, q_n]$. We also denote the usual cup product $a \smile b$ by
  $ab$ for all $a,b\in H^*(M)$.

  \medskip
  
  \proof[Proof of Theorem \ref{main}]
  
  This is an immediate consequence of the next lemma.
  
  \QED

  \begin{lemma}\label{l:main}
    Let $I= \{ 1\le i_1< i_2< \dots < i_k\le n\}$, and let $1\le i \le n$. Then	
    \begin{equation}
      \x{i_1}* \dots * \x{i_k} = \x{I} \text{ and }
      \x{i}*\x{i}=\1\ee{i}=q_i
    \end{equation}
  \end{lemma}
  
  \proof{} To prove the first equality we will proceed by induction. Assume 
  we have only two elements, say $\x{i},\x{j}$, with $i\ne j$. Then, by Proposition
  \ref{jump} and Remark \ref{rmk:jump} we have 
  $$
  \x{i}*\x{j}= x_{\{ij\}} + c \ \1  \otimes e^{B},
  $$
  where the coefficient $c$ is a rational number and $c_1(B)>0$.  
  
  From Corollary \ref{cor:prod} and the associativity of quantum
  multiplication we get
  \begin{align} 
    (\xp*\x{i})*\x{j}&=(\xc{\{i\}}*\x{j})\ee{i}\\
    &= \xc{\{ij\}} \ee{ij} + c\  \xp \otimes e^{B}.  
    \notag
  \end{align}
  By Proposition \ref{jump} the term $\xc{\{i\}}*\x{j}$ is of the form 
  $$
  \xc{\{i\}}\x{j} + \sum_{c_1(B') >0}{a_{B'}  \otimes e^{B'}} 
  $$
  where again $\deg(a_{B'})=\deg (\xc{\{i\}}) + \deg (\x{j})-2c_1(B')<2n$.
  Since $j\in\{i\}^c$, the term $\xc{\{i\}}\x{j}$ is zero. Thus we have 
  $$
   \sum_{c_1(B')>0}{a_{B'}  \otimes e^{B'}\ee{i}}=\xc{\{ij\}} \ee{ij} + c \ \xp \otimes
   e^B.  
  $$

  Then by comparing the degree of the coefficients in the previous equation,  the constant
  $c$ must vanish. 
  
  For the general case we will use the same argument. 
  Assume the result holds for $k$ different elements. 
  Let  $I'=\{i_{k+1} \}\cup I$. The quantum product  $\x{i_1}* \dots * \x{i_{k+1}}$
  is by the inductive hypothesis, the same as  $\x{I}*\x{i_{k+1}}$. This
  element can be written in terms of the basis as 
  $$
  \x{I}*\x{i_{k+1}}= \x{I'} + \sum_{c_1(B)>0,J\subset \mathcal{S}} a_{B,J}\  x_J \otimes e^B
  $$
  where $2|J|=\deg (x_J) = \deg (\x{I'}) -2d  \le \deg (\x{I'}) - 4$ and the coefficients
  $a_{B,J}$ are rational.

  As before, using quantum associativity and  Corollary \ref{cor:prod} we get
  \begin{align}\label{eq:d}
    (\xp*\x{I})*\x{i_{k+1}}&=(\xc{I}*\x{i_{k+1}})\ee{I}\\
    &=\xc{I'}\ee{I'}+ \sum_{c_1(B)>0,J\subset \mathcal{S} }{a_{B,J}\  \xc{J}  \otimes e^{A_J
    + B}}. 
    \notag
  \end{align}
  Here the degree satisfies
  \begin{equation}\label{eq:e}
  \deg (\xc{J}) = 2n - \deg (x_{I'}) + 2d \ge 2n - \deg (\x{I'}) +4=2(n-|I|+1).
  \end{equation}
  
  Now, the center term in  Equation (\ref{eq:d}) is written as
  
  $$
  (\xc{I}\x{i_{k+1}} + \sum_{c_1(B')>0, K\subset \mathcal{S} }{c_{B',K} \ x_K \otimes
  e^{B'} })\ee{I}, 
  $$
  where we have
  \begin{equation}\label{eq:f}
  \deg (x_K)\le \deg (\xc{I}) + \deg (\x{i_{k+1}}) - 4=2(n-|I|-1).
  \end{equation}
  Since $i_{i+1}\in I^c$, $\xc{I}\x{i_{k+1}}=0$. Finally we have the identity 
  $$
  \sum_{c_1(B')>0,K\subset \mathcal{S}}{ c_{B',K}  \ x_K\otimes  e^{B'+A_I} } = \xc{I'}
  \ee{I'}+ 
  \sum_{c_1(B)>0,J\subset \mathcal{S} }{a_{B,J} \ \xc{J}  \otimes e^{A_J + B}}.
  $$
  By Equations (\ref{eq:e}),(\ref{eq:f}),  the coefficients $a_{B,J}$ are zero. This
  proves the first part of the lemma.

  The second part is analogous, just write
  
  $$
  \x{i}*\x{i}= x_{i}\x{i} + c \ \1  \otimes e^B =c \ \1  \otimes e^B
  $$
  then multiplying by $\xp$
  $$
  (\xp*\x{i})*\x{i}=(\xc{\{i\}}*\x{i})\ee{i}=c \ \xp  \otimes e^B.
  $$
  Since $\xc{\{i\}}*\x{i}=\xp$, it follows that $c=1$ and $e^B=e^{A_i}$. 
  
  \QED


\begin{thebibliography}{999999999}
    
    %
  \bibitem{AuB} D. Austin and P. Braam, Morse--Bott theory and 
    equivariant cohomology, in {\it The Floer Memorial Volume},
    Progress in Mathematics {\bf 133}, Birkh\"auser (1995).
    %
    %
    %
    %
 

  \bibitem{Ha} A.Hattori, Symplectic manifolds with semifree Hamiltonian $S^1$ 
    actions,
    {\it Tokyo J. Math} {\bf 15} (1992) 281-296.
 
  \bibitem{K}    M. F. Kirwan, Cohomology of quotients in Symplectic and Algebraic Geometry.
          Mathematical Notes 31. {\it  Princeton University Press}, (1984).
    %
    %
    %
 
    %
    %
    %
    %
    %
    
  \bibitem{McT}  D. McDuff and S. Tolman, Topology of Hamiltonian circle 
    actions,
    {\it in preparation}.
    
    
  \bibitem{McS1}  D. McDuff and D. Salamon, {\it Introduction to
    Symplectic Topology}, 2nd edition (1998) OUP, Oxford, UK
    
  \bibitem{McS2}  D. McDuff and D. Salamon, {\it J-holomorphic curves and 
    Quantum Cohomology}, 2nd edition (2003) , AMS University Lecture Series, Vol {\bf 6}.
    
    
    
    %
    
    
    
    
    
  \bibitem{SZ}
    D. Salamon, and E. Zehnder,  Morse theory for periodic
    solutions of Hamiltonian systems and the Maslov index. {\it
      Communications in Pure and Applied Mathematics\/}, {\bf 45} (1992),
    1303--60.
    
    
    
  \bibitem{Sch} M. Schwarz, Equivalences for Morse homology, in {\it
    Geometry and Topology in Dynamics} ed M. Barge, K. Kuperberg,
    Contemporary Mathematics {\bf  246}, Amer. Math. Soc. (1999), 197--216.
  \bibitem{Sch2} M. Schwarz, {\it Morse Homology}, 1999, Birk\"auser Verlag.
     
    
  \bibitem{Sei}
    P. Seidel, $\pi_1$  of symplectic automorphism groups
    and invertibles in quantum cohomology rings, {\it Geom. and Funct.
      Anal.} {\bf 7} (1997), 1046 -1095.
    
    
  \bibitem{TW}
    S. Tolman and J. Weitsman, On semifree symplectic circle actions
    with isolated fixed points, {\it Topology}, {\bf 39} (2000), 299-309.
      
  \end{thebibliography}
\end{document}